\documentclass[10pt]{article}

\usepackage{amsmath,amscd,amsfonts,amsthm}

\newcommand{\shrinkmargins}[1]{
  \addtolength{\textheight}{#1\topmargin}
  \addtolength{\textheight}{#1\topmargin}
  \addtolength{\textwidth}{#1\oddsidemargin}
  \addtolength{\textwidth}{#1\evensidemargin}
  \addtolength{\topmargin}{-#1\topmargin}
  \addtolength{\oddsidemargin}{-#1\oddsidemargin}
  \addtolength{\evensidemargin}{-#1\evensidemargin} }

\shrinkmargins{.7}

\DeclareMathOperator{\Rea}{Re}
\DeclareMathOperator{\Ima}{Im}

\newcommand{\field}[1]{\mathbb{#1}}

\newcommand{\Q}{\field{Q}}

\newcommand{\Z}{\field{Z}}

\newcommand{\tensor} {\otimes}

\newcommand{\dual}[1]{#1^{\vee}}

\newcommand{\FF}{\mathcal{F}}
\newcommand{\notdiv}{\not | \,}
\newcommand{\beq}{\begin{displaymath}}
\newcommand{\eeq}{\end{displaymath}}
\newcommand{\beqn}{\begin{equation}}
\newcommand{\eeqn}{\end{equation}}

\theoremstyle{plain}
\newtheorem{thm}{Theorem}
\newtheorem{prop}[thm]{Proposition}
\newtheorem{cor}[thm]{Corollary}
\newtheorem{lem}[thm]{Lemma}

\theoremstyle{definition}

\theoremstyle{remark}

\title{On the error term in Duke's estimate for the average special value
of $L$-functions}   

\author{Jordan S. Ellenberg \\ Princeton University \\\texttt{ellenber@math.princeton.edu}}
\date{20 Apr 2005}

\begin{document}

\maketitle

\begin{abstract}  Let $\FF$ be an orthonormal basis for weight $2$
cusp forms of level $N$.  We show that various weighted averages of
special values $L(f \tensor \chi, 1)$ over $f \in \FF$ are equal to $4
\pi c + O(N^{-1 + \epsilon})$, where $c$ is an explicit nonzero constant. A
previous result of Duke gives an error term of $O(N^{-1/2}\log N)$.

\medskip

\end{abstract}

\section*{Introduction}

Let $N$ be a positive integer, and let $\FF$ be an basis for
$S_2(\Gamma_0(N))$ which is orthonormal for the Petersson inner
product.  Let $\chi$ be a Dirichlet character.

In \cite{duke:average}, Duke proves the estimate
\begin{equation}
\sum_{f \in \FF} a_1(f) L(f \tensor \chi, 1) = 4\pi + O(N^{-1/2} \log
N)
\label{eq:lsum}
\end{equation}
in case $N$ is prime and $\chi$ is unramified at $N$, using the
Petersson formula and the Weil bounds on Kloosterman sums. 

In this note, we will sharpen the error term in Duke's estimate to
$O(N^{-1+\epsilon})$.  At the same time, we observe that his techniques
generalize to arbitrary $N$ and $\chi$, and to the situation where
$a_1$ is replaced by an arbitrary $a_m$.

We have in mind an application to the problem of finding
all primitive solutions to the generalized Fermat equation
\begin{equation}
A^4 + B^2 = C^p
\label{eq:a4b2cp}
\end{equation}
In \cite{elle:a4b2cp}, we show how to associate to a solution of
\eqref{eq:a4b2cp} an elliptic curve over $\Q[i]$ with an isogeny to
its Galois conjugate and a non-surjective mod $p$ Galois
representation.  Such curves are parametrized by rational points on a
certain modular curve $X$; following Mazur's method, we can place
strong constraints on $X(\Q)$ by exhibiting a quotient of the Jacobian
of $X$ with Mordell-Weil rank $0$.  This problem, in turn, reduces via
the theorem of Kolyvagin and Logachev to proving the existence
of a new form $f$ on level $p^2$ or $2p^2$ such that the image of $f$
under a certain Hecke operator has an $L$-function with non-vanishing
special value.  We can then derive from Duke's estimate that
\eqref{eq:a4b2cp} has no solutions for $p > 2 \cdot 10^5$.  Using the
sharper estimate derived here, we find in \cite{elle:a4b2cp} that
\eqref{eq:a4b2cp} has no solutions for $p \geq 211$.

The author thanks Emmanuel Kowalski for useful discussions about the
topic of this paper, and is very grateful to Nathan Ng for finding an
error in an earlier version, and for suggesting several helpful
sharpenings of the bounds. 

\section*{Theorem statements}

In this section we state various versions of our estimate.  If $f$ is a
modular form, we always use $a_m(f)$ to denote the Fourier
coefficients of the $q$-expansion of $f$:
\beq
f = \sum_{m=0}^\infty a_m(f) q^m.
\eeq
As above, we denote by $\FF$ a Petersson-orthonormal basis for
$S_2(\Gamma_0(N))$.  

Write $(a_m, L_\chi)$ for the sum
\beq
\sum_{f \in \FF} a_m(f) L(f \tensor \chi, 1)
\eeq
and let $q$ be the conductor of $\chi$.

We obtain a rather complicated bound for $(a_m, L_\chi)$, which we
state below.

\begin{thm}  Suppose $N \geq 400$, $N \notdiv q$, and let $\sigma$ be a real number with $q^2/2\pi \leq \sigma \leq Nq/\log N$.  Then we can write
\beq
(a_m, L_\chi) = 4\pi\chi(m)e^{-2\pi m/ \sigma N \log N} - E^{(3)} +
E_3 - E_2 - E_1 + (a_m, B(\sigma N \log N))
\eeq
where
\begin{itemize}
\item $|(a_m,B(\sigma N \log N)| \leq 30 (400/399)^3 \exp(2\pi) q^2 m^{3/2}
N^{-1/2} d(N) N^{-2 \pi \sigma / q^2}$;
\item $|E_1| \leq (16/3)\pi^3 m^{3/2} \sigma \log N e^{-N/2\pi m \sigma
\log N}$;
\item $|E_2| \leq (8/9)\pi^5\zeta^2(7/2)m^{5/2} \sigma^2 N^{-3/2}\log^2 N$;
\item $|E_3| \leq (8/3) \zeta^2(3/2) \pi^3 \sigma m^{3/2} N^{-1/2} \log N
d(N) e^{-N/ 2 \pi m \sigma \log N}$;
\item
$|E^{(3)}| \leq 16\pi^3 m \sum_{c>0,N|c}
\min[\frac{2}{\pi} \phi(q) c^{-1} \log c , \frac{1}{6}\sigma N \log N m^{1/2}c^{-3/2}d(c)].$

\end{itemize}
\label{th:main}
\end{thm}

\begin{proof}  Immediate from
Propositions~\ref{pr:ab},\ref{pr:e1},\ref{pr:e2},\ref{pr:e3},\ref{pr:eu3}.
\end{proof}

If $q,m$ are considered as constants, the bound above simplifies considerably.

\begin{cor}
\beq
(a_m,L_\chi) = 4\pi\chi(m)e^{-2\pi m/ \sigma N \log N} + O(N^{-1+\epsilon})
\eeq
where the implied constants depend only on $m,q,$ and $\epsilon$.
\end{cor}

\begin{proof} The only thing to check is that the bound on $|E^{(3)}|$
  is of order at most $N^{-1+\epsilon}$; one checks this by fixing
  some cutoff $X$, say $X=N^3$, and observing that both $\sum_{0<c<X,N|c} c^{-1} \log c$ and
  $N \log N \sum_{c>X,N|c} c^{-3/2}d(c)$ are $O(N^{-1+\epsilon})$.
\end{proof}

The ``true behavior'' of $(a_m,L_\chi)$ is less clear.  One might
for instance ask: what is the true asymptotic behavior of
$(a_m,L_\chi) - 4\pi\chi(m)$ as $N$ grows with $m,q$ held fixed?  More
generally, what is the shape of the region in $m,q,N$-space for which
$(a_m,L_\chi)$ is close to $4\pi\chi(m)$?  One might, for instance,
define $f_\delta(N)$ to be the smallest integer such that
$|(a_m,L_\chi) - 4 \pi\chi(m)| \leq \delta$ for all $m \leq f(N)$.
Duke's approach shows that $f_\delta(N) \gg N^{1/2}$, whereas the
present results show that $f_\delta(N) \gg N^{3/5}$.  (Remark: further
expansion of the Bessel function in Taylor series will give
$f_\delta(N) \gg N^{1-\epsilon}$, with a constant depending on
$q,\epsilon$.)  Similarly, one could try to optimize the dependence on
$q$ in order to get a result that applied when $q$ is large compared to $N$.



\section*{Proof of the main result}

We begin by recalling the Petersson trace formula.  

\begin{lem}[Petersson trace formula]  Let $m,n$ be positive integers,
and let $\FF$ be an orthonormal basis for $S_2(\Gamma_0(N))$.

Then
\begin{equation}
\frac{1}{4\pi\sqrt{mn}}\sum_{f \in \FF} a_m(f)a_n(f) =  
\delta_{mn} - 2\pi \sum_{\stackrel{c > 0} {c = 0 \pmod N}}
c^{-1}S(m,n;c)J_1(4\pi\sqrt{mn}/c)
\label{eq:petersson}
\end{equation}
where $S(m,n;c)$ is the Kloosterman sum for $\Gamma_0(N)$, and $J_1$
is the $J$-Bessel function. 
\label{le:petersson}
\end{lem}

\begin{proof} See \cite[Th.\ 3.6]{iwan:caf}. \end{proof}

We can and do assume that $\FF$ consists of eigenforms for $T_p$ for
all $p \notdiv N$, and for $w_N$.

The Petersson product on $S_2(\Gamma_0(N))$ induces an inner product
on the dual space $\dual{S_2(\Gamma_0(N))}$.  With respect to this product, the
left-hand side of \eqref{eq:petersson} is
$\frac{1}{4\pi\sqrt{mn}}(a_m,a_n).$

Lemma~\ref{le:petersson} immediately gives a bound on the size of $(a_m,a_n)$.

\begin{lem} We have the bound
\beq
|(a_m,a_n) - 4\pi\sqrt{mn}\delta_{mn}| \leq
8 \zeta^2(3/2) \pi^2 (m,n)^{1/2} mn N^{-3/2} d(N).
\eeq
\label{le:amanbound}
\end{lem}

\begin{proof}
Applying the Weil bound
\beq
|S(m,n;c)| \leq (m,n,c)^{1/2}d(c)c^{1/2}
\eeq
and the fact that $|J_1(x)| \leq x/2$ yields
\begin{eqnarray*}
|4\pi\sqrt{mn}\sum_{\stackrel{c > 0} {c = 0 \pmod N}}
c^{-1}S(m,n;c)J_1(4\pi\sqrt{mn}/c)|  & \leq &
4 \pi \sqrt{mn} \sum_{\stackrel{c > 0} {c = 0 \pmod N}}
c^{-1/2} d(c) (m,n)^{1/2} (2 \pi \sqrt{mn} / c) \\
& = & 8 \pi^2 (m,n)^{1/2} mn \sum_{\stackrel{c > 0} {c = 0 \pmod N}}
c^{-3/2} d(c).
\end{eqnarray*}
Now the sum over $c$ is equal to
\beq
\sum_{b > 0} (Nb)^{-3/2} d(Nb)
\eeq
which is bounded above by
\beq
N^{-3/2} d(N) \sum_{b>0} b^{-3/2} d(b) = \zeta^2(3/2) N^{-3/2} d(N).
\eeq
This yields the desired result.
\end{proof}

Let $L_\chi$ be the element of $\dual{S_2(\Gamma_0(N))}$ which sends
each cusp form $f$ to the special value $L(f \tensor \chi, 1)$.  Then
the value to be estimated is precisely $(L_\chi,a_m)$.  In order to
estimate this product via the Petersson formula, it is necessary to
approximate $L_\chi$ as a sum of Fourier coefficients.  We accomplish
this via the standard approximation to $L_\chi(f)$ by a rapidly
converging series~\cite{rohr:cycl}.

We define a linear functional $A(x)$ on
$S_2(\Gamma_0(N))$ by the rule
\beq
A(x)(f) = \sum_{n \geq 1} \chi(n) a_n(f) n^{-1} e^{-2\pi n / x}.
\eeq
Then $A$ is a good approximation to the functional $L_{\chi}$ when $x$
becomes large.  Let $B(x) = A(x) - L_\chi$.  Let $M$ be an integer
such that $f \tensor \chi$ is a cuspform on $\Gamma_1(M)$ for all $f \in \FF$.

By the functional equation for $L(f \tensor \chi,s)$, we have
\beq
B(x)(f) = \sum_{n \geq 1} a_n(w_M(f \tensor \chi)) n^{-1}
e^{-2\pi nx/M}.
\eeq

When $x$ is on the order of $N \log N$, then $B(x)$ is a short sum,
and we want to show it is negligible.  The only difficulty is bounding
the Fourier coefficients of $w_M(f \tensor \chi)$.  This is difficult
only in case the conductor of $\chi$ has common factors with $N$, in
which case $f \tensor \chi$ is not necessarily an eigenform for any
$W$-operator, even when $f$ is a new form (see \cite{atki:twistw}.)

A crude bound will be enough for us.  We define an ``average cuspform'' 
\beq
g = \sum_{f \in \FF} a_m(f) (f \tensor \chi).
\eeq
Then
\beq
a_n(g) = \chi(n) (a_m,a_n) 
\eeq
and it follows from Lemma~\ref{le:amanbound} that
\beq
|a_n(g)| \leq (8 \zeta^2(3/2) \pi^2 m^{3/2} N^{-3/2} d(N))n
\eeq
for all $n \neq m$, while
\beq
|a_m(g)| \leq 4\pi\sqrt{mn} + (8 \zeta^2(3/2) \pi^2 m^{3/2} N^{-3/2} d(N))n
\eeq
when $m=n$.

We have that
\beq
(a_m,B(x)) = \sum_{f \in \FF} a_m(f) \sum_{n > 0} a_n(w_M (f
\tensor \chi)) n^{-1} e^{-2\pi nx/M} =
\sum_{n > 0} a_n(w_M g) n^{-1} e^{-2\pi nx/M},
\eeq
so it remains to bound the Fourier coefficients of the single form
$w_M g$.  Write $c$ for the constant $8 \zeta^2(3/2) \pi^2 m^{3/2}
N^{-3/2} d(N)$.


If $\tau$ is a point in the upper half plane, we have
\begin{eqnarray*}
|g(\tau)| \leq \sum_{n > 0} |a_n e^{2\pi i n \tau}| & = & \sum_{n > 0}
|a_n| 
\exp(-2\pi\Ima(n \tau)) \\
& \leq & 
\sum_{n > 0} cn \exp(-2\pi\Ima(n \tau))
+ 4\pi m \exp(-2\pi\Ima(m \tau)) \\
& \leq & c (2 \pi \Ima(\tau))^{-2} + 4 \pi m.
\label{eq:gtau}
\end{eqnarray*}

Choose a positive real constant $\alpha$. 
The Fourier coefficient $a_n(w_M g)$ can be expressed as 
\begin{equation}
\int^1_0 w_M g (\alpha i + t) \exp(-2\pi i n (\alpha i + t)) dt 
= \int^1_0 M^{-1} (\alpha i + t)^{-2} g(-1/M(\alpha i + t)) 
\exp(-2\pi i n (\alpha i + t)) dt.
\label{eq:anint}
\end{equation}
Now $\Ima((-1/M(\alpha i + t))) = M^{-1}\alpha |\alpha i + t|^{-2}$.
So it follows from \eqref{eq:gtau} that
\begin{eqnarray*}
|a_n(w_M g)| & \leq &
\int^1_0 M^{-1}  |\alpha i + t|^{-2}  [c (2\pi)^{-2} M^2 \alpha^{-2}
|\alpha i + t|^4 + 4 \pi m] \exp(2 \pi n \alpha) dt \\
& = & c M (2\pi)^{-2} \exp(2 \pi i n \alpha) \alpha^{-2} \int^1_0
|\alpha i + t|^{2} dt + 4 \pi m M^{-1} \exp(2 \pi n \alpha)
\int^1_0 |\alpha i + t|^{-2} dt \\
& \leq & c M (2\pi)^{-2} \exp(2 \pi n \alpha) \alpha^{-2} (\alpha^2
+ 1) + 4 \pi m M^{-1} \exp(2 \pi n \alpha) \alpha^{-2}.
\end{eqnarray*}

Now setting $\alpha = 1/n$ yields
\beq
|a_n(w_M g)| \leq c M (2\pi)^{-2} \exp(2 \pi)(1 + n^2) + 4 \pi
\exp(2\pi) m M^{-1} n^2.
\eeq

We now use the very rough bound $1+n^2 \leq n^2(n+1)$ to obtain
\begin{eqnarray*}
|(a_m, B(x))| & = & |\sum_{n > 0} a_n(w_M g) n^{-1} e^{-2\pi
nx/M}| \\
& \leq &
[c M (2\pi)^{-2} \exp(2\pi) + 4 \pi m M^{-1} \exp(2\pi)]
\sum_{n > 0} n(n+1)e^{-2\pi n x / M} \\
& = & \exp(2\pi)(c M (2\pi)^{-2} + 4 \pi m M^{-1})(2 \exp(-2\pi x /
M))(1-\exp(-2\pi x / M))^{-3}.
\end{eqnarray*}

Now $M$ can be taken to be $q^2 N$ where $q$ is the conductor of
$\chi$.  Let $\sigma$ be a constant to be fixed later, and set $x =
\sigma N \log N$.  Finally, suppose $N > 400$ and suppose $\sigma >
q^2/2\pi$.  First of all, we observe that under the hypothesis on $N$,
\begin{eqnarray*}
c M (2\pi)^{-2} + 4 \pi m M^{-1}
& = & 2 \zeta^2(3/2) q^2  m^{3/2} N^{-1/2} d(N) + 4 \pi m q^{-2} N^{-1} \\
& \leq & 
15 q^2 m^{3/2} N^{-1/2} d(N).
\end{eqnarray*}
Also,
\beq
1-\exp(-2\pi x / M) = 1 - \exp(-2 \pi \sigma \log N / q^2) 
\leq 1 - 400^{-2 \pi \sigma/q^2} \leq 400/399.
\eeq
So, in all, we have proved the following.

\begin{prop} Suppose $N \geq 400$ and $\sigma > q^2/2\pi$.  Then
\beq
|(a_m,B(\sigma N \log N)| \leq 30 (400/399)^3 \exp(2\pi) q^2 m^{3/2}
N^{-1/2} d(N) N^{-2 \pi \sigma / q^2}.
\eeq
\label{pr:ab}
\end{prop}

In other words, we have shown that the error in approximating
$(a_m,L_\chi)$ by $(a_m,A(x))$ is bounded by a function decreasing
quickly in $N$, if $x$ is chosen on the order of $q^2 N \log N$.

We now turn to the analysis of $(a_m, A(\sigma N \log N))$.

\medskip

First of all, we have
\beq
(a_m,A(\sigma N \log N))
= \sum_{f \in \FF} a_m(f) \sum_{n > 0} \chi(n) a_n(f) n^{-1} e^{-2\pi
n / \sigma N \log N}
= \sum_{n > 0} \chi(n) (a_m,a_n) n^{-1} e^{-2\pi
n / \sigma N \log N}
\eeq
which, by Lemma~\ref{le:petersson}, equals
\beq
4\pi \chi(m) e^{-2 \pi m / \sigma N \log N}
- 8 \pi^2 \sqrt{m} \sum_{n > 0} \chi(n) n^{-1/2} e^{-2\pi
n / \sigma N \log N} \sum_{\stackrel{c > 0} {c = 0 \pmod N}}
c^{-1}S(m,n;c)J_1(4\pi\sqrt{mn}/c).
\eeq
We split the latter sum into two ranges; write
\beq
E^{(1)} = 8 \pi^2 \sqrt{m} \sum_{n > 0} \chi(n) n^{-1/2} e^{-2\pi
n / \sigma N \log N} \sum_{\stackrel{c > 2\pi\sqrt{mn}} {c = 0 \pmod N}}
c^{-1}S(m,n;c)J_1(4\pi\sqrt{mn}/c)
\eeq
and
\beq
E_1 = 8 \pi^2 \sqrt{m} \sum_{n > 0} \chi(n) n^{-1/2} e^{-2\pi
n / \sigma N \log N} \sum_{\stackrel{0 < c \leq 2\pi\sqrt{mn}} {c = 0 \pmod N}}
c^{-1}S(m,n;c)J_1(4\pi\sqrt{mn}/c).
\eeq
We claim $E_1$ decreases quickly with $N$.  First, recall that
$|J_1(a)| \leq \min(1,a/2)$ for all real $a$.  So
\beq
|E_1| \leq 8 \pi^2 \sqrt{m} \sum_{n > 0} n^{-1/2} e^{-2\pi
n / \sigma N \log N}\sum_{0 < Nb \leq 2\pi\sqrt{mn}}
(Nb)^{-1} S(m,n;Nb).
\eeq
Note that the inner sum in $|E_1|$ has nonzero terms only when $n >
(N/2\pi\sqrt{m})^2$.  In this range, the exponential decay takes
over.  
We observe that $|S(m,n;Nb)| \leq m^{1/2}(Nb)^{1/2}d(Nb) <
2\sqrt{m}Nb$, so we can bound $E_1$ by  
\begin{eqnarray*}
|E_1| & \leq & 8 \pi^2 \sqrt{m} \sum_{n > (N/2\pi\sqrt{m})^2} n^{-1/2} e^{-2 \pi n /
\sigma N \log N} \sum_{0 < Nb \leq 2\pi\sqrt{mn}} 2\sqrt{m} \\
& \leq &
8 \pi^2 \sqrt{m} \sum_{n > (N/2\pi\sqrt{m})^2} n^{-1/2} e^{-2 \pi n /
\sigma N \log N} (2\sqrt{m}) (2 \pi \sqrt{mn}/N) \\
& = & 
32 \pi^3 N^{-1} m^{3/2} \sum_{n > (N/2\pi\sqrt{m})^2} e^{-2 \pi n /
\sigma N \log N} \\
& \leq & 32 \pi^3 N^{-1} m^{3/2} e^{-N / 2 \pi m \sigma \log N}(1-e^{-2
\pi / \sigma N \log N})^{-1}.
\end{eqnarray*}

We now simplify this bound under assumptions on $N$ and $\sigma$.
\begin{prop} Suppose $N \geq 400$ and $\sigma > q^2/2\pi$.  Then
\beq
|E_1| \leq (16/3)\pi^3 m^{3/2} \sigma \log N e^{-N/2\pi m \sigma \log N}.
\eeq
\label{pr:e1}
\end{prop}

\begin{proof} This amounts to the observation that $\sigma N \log N
\geq 300$, from which it follows that
\beq
(1-e^{-2 \pi / \sigma N \log N})^{-1} \leq (1/6)\sigma N \log N.
\eeq
\end{proof}


\medskip

We now consider the sum $E^{(1)}$ over the range where $n$ is small
compared to $c$.  In this range, we use the Taylor approximation 
\begin{equation}
|J_1(a) - a/2| \leq (1/16)a^3.
\label{eq:besselt}
\end{equation}
So we can write $E^{(1)} = E^{(2)} + E_2$, where
\beq
E^{(2)} = 8 \pi^2 \sqrt{m} \sum_{n > 0} \chi(n) n^{-1/2} e^{-2\pi
n / \sigma N \log N} \sum_{\stackrel{c > 2\pi\sqrt{mn}} {c = 0 \pmod N}}
c^{-1}S(m,n;c)(2\pi\sqrt{mn}/c).
\eeq
We claim $E_2$ decreases with $N$.  For we have by \eqref{eq:besselt} that
\begin{eqnarray*}
|E_2| & \leq & 8 \pi^2 \sqrt{m} \sum_{n > 0} n^{-1/2} e^{-2\pi
n / \sigma N \log N} \sum_{\stackrel{c > 2\pi\sqrt{mn}} {c = 0 \pmod N}}
c^{-1}S(m,n;c)(1/16)(4\pi\sqrt{mn}/c)^3 \\
& = & 32 \pi^5 m^2 \sum_{n > 0}\sum_{\stackrel{c > 2\pi\sqrt{mn}} {c = 0
\pmod N}} n e^{-2\pi n / \sigma N \log N} \sum_{\stackrel{c >
2\pi\sqrt{mn}} {c = 0 \pmod N}} c^{-4} S(m,n;c).
\end{eqnarray*}
We now use the Weil bound $S(m,n;c) \leq m^{1/2} c^{1/2} d(c)$ to get 
\begin{eqnarray*}
|E_2| & \leq & 32 \pi^5 m^{5/2} \sum_{n > 0}\sum_{\stackrel{c > 2\pi\sqrt{mn}}
{c = 0 \pmod N}}  
n e^{-2\pi n / \sigma N \log N}
c^{-7/2} d(c)
\\
& \leq & 
32 \pi^5 m^{5/2} \sum_{n > 0} \sum_{b > 0} n e^{-2\pi n / \sigma N
\log N} N^{-7/2} d(N) b^{-7/2} d(b) \\
& \leq &
32 \pi^5 m^{5/2} N^{-7/2} d(N) \zeta^2(7/2) \sum_{n > 0} n e^{-2\pi n
/ \sigma N \log N}
\end{eqnarray*}
So we can write
\beq
|E_2| \leq 32 \pi^5 \sqrt{3} \zeta(3) m^{5/2} N^{-7/2} e^{-2\pi / \sigma N \log N} (1-
e^{-2\pi / \sigma N \log N})^{-2}.
\eeq

\begin{prop} Suppose $N > 400$ and $\sigma > q^2/2\pi$.  Then
\beq
|E_2| \leq (8/9)\pi^5\zeta^2(7/2)m^{5/2} \sigma^2 N^{-3/2}\log^2 N.
\eeq
\label{pr:e2}
\end{prop}

\begin{proof} Another use of the bound $(1- e^{-2\pi / \sigma N \log
N})^{-1} \leq (1/6)\sigma N \log N$.
\end{proof}

\medskip


 

We now come to $E^{(2)}$, which is the main term of the error
\beq
|(a_m, L_\chi) - 4\pi\chi(m)e^{-2\pi m/\sigma N \log N}|.
\eeq
Recall from above that 
\beq
E^{(2)} = 16 \pi^3 m \sum_{n > 0} \sum_{\stackrel{c > 2\pi\sqrt{mn}}
{c = 0 \pmod
N}} \chi(n) e^{-2\pi
n / \sigma N \log N} c^{-2}S(m,n;c).
\eeq

Applying the Weil bound to $S(m,n;c)$ yields the estimate $E^{(2)} =
O(N^{-1/2} \log N)$ which appears in \cite{duke:average}.  We want to
exploit cancellation between the Kloosterman sums in order to improve
Duke's bound on $E^{(2)}$.

For simplicity, we carry this out under assumptions on the size of $N$
and $\sigma$.  For the remainder of this section, assume that
\begin{itemize}
\item $N \geq 400$;
\item $q^2/2\pi \leq \sigma \leq Nq/\log N$.
\end{itemize}

Recall that under these hypotheses
\beq
\sigma N \log N \geq (1/2\pi) 400 \log 400 > 300.
\eeq

First of all, we will need a simple bound on the modulus of $1-e^z$.

\begin{lem} Let $z$ be a complex number with $|\Ima z| \leq \pi$ and $-2\pi/30
\leq \Rea z \leq 0$.  Then 
\beq
(1/2)|z| \leq |1-e^z| \leq |z|.
\eeq
\label{le:ezbound}
\end{lem}

\begin{proof}  The extrema of $|1-e^z|/|z|$ lie on the boundary of the
rectangular region under consideration; now a consideration of the
derivatives of $|1-e^z|/|z|$ on each of the four edges of the region
shows that the extrema are at the corners.  Computation of the values
of $|1-e^z|/|z|$ gives the result.
\end{proof}

Write
\beq
E^{(3)} = 16 \pi^3 m \sum_{n > 0} \sum_{\stackrel{c > 0}
{c = 0 \pmod
N}} \chi(n) e^{-2\pi
n / \sigma N \log N} c^{-2}S(m,n;c)
\eeq
and
\beq
E_3 = 16 \pi^3 m \sum_{n > 0} \sum_{\stackrel{c \leq 2\pi\sqrt{mn}}
{c = 0 \pmod
N}} \chi(n) e^{-2\pi
n / \sigma N \log N} c^{-2}S(m,n;c).
\eeq

So $E^{(2)} = E^{(3)} - E_3$.

The sum $E_3$, like $E_1$,  is supported in the region where exponential
decay dominates.  To be precise, the inner sum in $E_3$ has nonzero terms only when $n \geq (c/2\pi\sqrt{m})^2 \geq
N^2/4\pi^2 m$.  It follows that
\begin{eqnarray*}
|E_3| & \leq & 16 \pi^3 m \sum_{n > N^2/4\pi^2 m} \sum_{\stackrel{c > 0}{c = 0 \pmod
N}}
e^{-2\pi n / \sigma N \log N} m^{1/2} c^{-3/2} d(c) \\
& \leq & 16 \zeta^2(3/2) \pi^3 m^{3/2} (N^{-3/2} d(N)) e^{-N/ 2 \pi m \sigma \log N}
(1-e^{-2\pi/\sigma N \log N})^{-1}.
\end{eqnarray*}

Using the lower bounds on $N$ and $\sigma$, we obtain
\begin{prop}  Suppose $N > 400$ and $\sigma > q^2/2\pi$.  Then
\beq
|E_3| \leq (8/3) \zeta^2(3/2) \pi^3 \sigma m^{3/2} N^{-1/2} \log N d(N) e^{-N/ 2 \pi m
\sigma \log N}.
\eeq
\label{pr:e3}
\end{prop}

It now remains only to bound the main term
\beq
E^{(3)} = 16 \pi^3 m \sum_{n > 0} \sum_{\stackrel{c > 0}
{c = 0 \pmod
N}} \chi(n) e^{-2\pi
n / \sigma N \log N} c^{-2}S(m,n;c)
\eeq
We can write
\begin{equation}
E^{(3)} = 16 \pi^3 m \sum_{\stackrel{c > 0} {c = 0 \pmod N}} c^{-2}
S(c)
\label{eq:eu3}
\end{equation}
where 
\begin{eqnarray*}
S(c) & = & \sum_{n > 0} \chi(n) e^{-2\pi n / \sigma N \log N} S(m,n;c)
\\
& = & \sum_{x \in (\Z/c\Z)^*} \sum_{n > 0} \chi(n) e^{-2 \pi n / \sigma N
\log N} e\left(\frac{mx + ny}{c}\right)
\end{eqnarray*}
where $e(z) = e^{2 \pi i z}$ and $y \in (Z/c\Z)^*$ is the
multiplicative inverse of $x$.

For ease of notation, write $A = \sigma N \log N$, and for each
integer $y$ write $\epsilon_y
= 2 \pi (-1/A + yi/c)$.  Then
\begin{eqnarray*}
|S(c)| & \leq & \sum_{x \in (\Z/c\Z)^*} | \sum_{n > 0} \chi(n) e^{-2 \pi n/A} e\left(\frac{ny}{c}\right)| \\
& = & \sum_{x \in (\Z/c\Z)^*} | \sum_{\alpha = 1}^q \chi(\alpha)
e^{-2 \pi \alpha / A} e\left(\frac{\alpha y}{c}\right) \sum_{\nu \geq 0}
e^{2 \pi q\nu/A} e\left(\frac{q\nu y}{c}\right)| \\
& = & \sum_{x \in (\Z/c\Z)^*} | \sum_{\alpha=1}^q \chi(\alpha)
e^{-2 \pi \alpha / A} e\left(\frac{\alpha y}{c}\right) (1 - e^{2\pi q(-1/A +
iy/c)})^{-1} |\\
& = & \sum_{y \in (\Z/c\Z)^*} | (1-e^{q\epsilon_y})^{-1} \sum_{\alpha=1}^q
\chi(\alpha) e^{\alpha \epsilon_y}| \\
& \leq & \sum_{y \in (\Z/c\Z)^*} |(1-e^{q\epsilon_y})|^{-1} |\sum_{\alpha=1}^q
\chi(\alpha) e^{\alpha \epsilon_y}|.
\end{eqnarray*}

We have the trivial bound $|\sum_{\alpha=1}^q
\chi(\alpha) e^{\alpha \epsilon_y}| \leq \phi(q)$.  (This bound can be
sharpened to $O(\sqrt{q} \log q)$ if one wishes to improve the
dependence on $q$.) We now estimate $|\sum_y
(1-e^{q\epsilon_y})^{-1}|$.  For each $y$, let $f(y)$ be the unique
integer congruent to $qy$ modulo $c$ with $|f(y)| \leq c/2$.   By our
assumption that $N \notdiv q$, we have $f(y) \neq 0$.  Then by
Lemma~\ref{le:ezbound} one has
\beq
|(1-e^{q\epsilon_y})^{-1}| < \frac{c}{\pi |f(y)|}.
\eeq
Now the values of $|f(y)|$ range over the integers $a$ between $1$ and
$c/2$ such that $(a,c) = (q,c)$, each of which arises from at most $2(q,c)$
values of $y$.  So we have
\beq
|\sum_{y \in (\Z/c\Z)^*} (1-e^{q\epsilon_y})^{-1}|
\leq
\frac{2(q,c)c}{\pi}\left[\frac{1}{(q,c)} + \frac{1}{2(q,c)} + \ldots +
  \frac{1}{r(q,c)}\right]
=
(2c/\pi)[1 + 1/2 + \ldots + 1/r]
\eeq
where $r$ is the largest integer such that $r(q,c) \leq c/2$.  The value of $(2c/\pi)[1 + \ldots + 1/r]$  is largest when $(q,c) = 1$; in that case it is
bounded above by
\beq
(2c/\pi)[\log(c/2) + \gamma + 2/c].
\eeq
Since $c > 400$, the above expression is bounded by
$(2/\pi)c \log c$.  So, in all, one has
\begin{equation}
|S(c)| < (2/\pi) \phi(q) c \log c .
\label{eq:smallrange}
\end{equation}

We observe as well that, from the Weil bound, we have
\beq
|S(c)| \leq  \sum_{n > 0} e^{-2\pi n / A} m^{1/2} c^{1/2} d(c) \leq
m^{1/2} c^{1/2} d(c) (1 - e^{-2\pi/A})^{-1}.
\eeq

Recall from the proof of Proposition~\ref{pr:e1} that $(1 - e^{-2\pi/A})^{-1}
\leq (1/6)A$ under our conditions on $N$ and $\sigma$.  So 
\begin{equation}
|S(c)| \leq (1/6)Am^{1/2}c^{1/2}d(c).
\label{eq:largerange}
\end{equation}

In particular, we immmediately have the following proposition:
 
\begin{prop}  Suppose $N \geq 400$, $N \notdiv q$, and $\sigma >
  q^2/2\pi$.  Then  
\beq
|E^{(3)}| \leq 16\pi^3 m \sum_{\stackrel{c > 0} {c = 0 \pmod N}}
\min[\frac{2}{\pi} \phi(q) c^{-1} \log c , \frac{1}{6}\sigma N \log N m^{1/2}c^{-3/2}d(c)].
\eeq
\label{pr:eu3}
\end{prop}

This completes the proof of Theorem~\ref{th:main}.

\end{document}